
\input amstex
\input amsppt.sty
\magnification=\magstep1
\hsize = 6.5 truein
\vsize = 9 truein
\NoBlackBoxes
\NoRunningHeads\refstyle{B}\widestnumber\key{KhST}
\TagsAsMath
\TagsOnRight
\catcode`\@=11
\redefine\logo@{}
\catcode`\@=13

\def\b{\beta}
\def\la{\lambda}\def\CC{\bold C}\def\ZZ{\bold Z}
\def\deg{\operatorname{deg}}\def\res{\operatorname{res}}\def\a{\alpha}
\def\b{\beta}

\topmatter      \title\nofrills
Compatible Poisson-Lie structures on the loop group of $SL_{2}$.
    \endtitle      \author  B. Enriquez and V. Rubtsov  \endauthor
     \address  \newline
B.E., V.R.:
Centre de Math\'{e}matiques, Ecole Polytechnique, URA 169 du CNRS,
\newline
F--91128 Palaiseau, France \newline
V.R.(permanent address): ITEP, B. Cheremushkinskaya, 117259,
Moscou, Russie
\endaddress

 \abstract
We define a 1-parameter family of $r$-matrices on the loop algebra of
$sl_{2}$, defining compatible Poisson structures on the associated
loop group, which degenerate into the rational and trigonometric
structures, and study the Manin triples associated to them.
\endabstract

\endtopmatter

\subheading{Introduction} The concept of bi-(or multi-)Hamiltonian
structures on manifolds is known, since the works of Magri and
Gelfand-Dorfman (\cite{Ma}, \cite{GD}), to play an important role in the
study of classical integrable systems. Up to now, it still remains
unclear whether this notion has any reasonable quantum analogue. Among
long-known examples of such structures are the KdV hierarchy phase
spaces. In that case, higher compatible Poisson structures are nonlocal,
and are expected to be quantized using (nonlocal) vertex operator
algebras. This program is very far from its goals at the moment.

An other example of bi-Hamiltonian manifolds comes from the theory of
Poisson-Lie groups. It is an elementary fact that the rational and
trigonometric structures (according to the terminology of \cite{BD}) on
loop groups, are compatible. Accordingly, it is a natural question to
find higher compatible Poisson structures on these groups.
In this paper, we answer this question in the $sl_{2}$ case. We perform
an elementary manipulation
on the rational $r$-matrix, that after the operation of shifting of the
spectral parameter generates the rational and trigonometric
$r$-matrices. This shows the compatibility of these Poisson
structures. We then study the associated Manin triples. They are of a
slightly different kind from those studied in \cite{D}.
It would be interesting to achieve explicitly the quantization of the
structures found here.

A possible connection between the two problems mentioned above could be the
remark, that the Poisson brackets of the KdV monodromy operators in the
first (resp. second) Hamiltonian structures are respectively given by
the rational and trigonometric $r$-matrices (\cite{FT}).

\subheading{1. Compatible $r$-matrices for $sl_{2}(\CC((\la^{-1})))$}

Let $e,f,h$ be the Chevalley basis of $a=sl_{2}(\CC)$, $t=h\otimes
h+2(e\otimes f +f\otimes e)$ the split Casimir element of
$sl_{2}(\CC)$.
The rational and trigonometric
$r$-matrices of $a(\CC((\la)))$ are
$$
r_1(\la,\mu)={t\over{\la-\mu}}\quad \text{and}\quad
r_2(\la,\mu)={{\la+\mu}\over{\la-\mu}}t+2(e\otimes f +f\otimes e)
$$
They can be considered as elements of $Usl_{2}(\CC)^{\otimes 2}\otimes
\CC(\la,\mu)$, and
satisfy the Yang-Baxter equation with parameters, that is
$$
[\![r,r]\!](\la,\mu,\nu):=
[r^{12}(\la,\mu), r^{13}(\la,\nu)]+\text{cyclic permutations}=0,
$$
in the algebra
$Usl_2(\CC)^{\otimes 3}\otimes
\CC(\la,\mu,\nu)$. Let us perform the following transformation on $r_1$:
we replace $\la$ and $\mu$ by their inverses and apply the affine Weyl
group
transformation $\la^ke\to\la^{k+1}e$, $\la^kh\to\la^{k}h$, $\la^kf\to
\la^{k-1}f$. Up to sign the resulting $r$-matrix is
$$
r_3(\la,\mu)={{\la\mu}\over{\la-\mu}}h\otimes h+2{\la^2\over{\la-\mu}}
e\otimes f+2{\mu^2\over{\la-\mu}}f\otimes e=
{{\la\mu}\over{\la-\mu}}t+2\la e\otimes f -2\mu f\otimes e
\eqno(1)
$$
Adding a quantity $E$ to $\la$ and $\mu$ transforms $r_3$ into $r_3+Er_2+E^2
r_1$, so $[\![r_3+Er_2+E^2r_1,r_3+Er_2+E^2r_1]\!]=0$. From this follows
that $[\![r_1,r_2]\!]=[\![r_2,r_3]\!]=0$ and also that
$[\![r_1,r_3]\!]=0$ (since we know that $[\![r_2,r_2]\!]=0$; we set
$[\![a,b]\!]=[a^{12}(\la,\mu),b^{13}(\la,\nu)]+[b^{12}(\la,\mu),
a^{13}(\la,\nu)]+$c.p.). Hence, the
Poisson structures defined on
$SL_2(\CC((\la)))$ by $r_1$, $r_2$ and $r_3$ are compatible (in the
sense of \cite{Ma}, \cite{GD}). Let $\tau^{(0)}_{E}$ be the automorphism of
$\CC((\la^{-1}))$ defined by
$\tau_{E}(\la^{k})=(\la+E)^{k}=\la^{k}+kE\la^{k-1}+...$ and let us call
$\tau_{E}$ the automorphism induced on $SL_{2}(\CC((\la^{-1})))$.
We have:

\proclaim{Proposition} For $a_{i}$ complex numbers ($i=1,2,3$), the
multiplicative bivectors defined by $\sum_{i=1}^{3}a_{i}r_{i}$ define
structures of
Poisson-Lie group on $SL_{2}(\CC((\la^{-1})))$. Up to multiplication and
the action of the automorphisms $\tau_{E}$, this family is reduced to a
1-parameter family, with special points the rational and trigonometric
structures.
\endproclaim

\subheading{2. Modified Yang-Baxter equation}

Let us remark that solutions to the Yang-Baxter equation with parameters
usually correspond to solutions of the modified Yang-Baxter equation,
after interpretation of
$r$-matrices
as elements of $[sl_2(\CC((\la^{-1})))\otimes sl_2(\CC((\mu^{-1})))]^-$
(${}^-$ means
direct product of homogeneous components, where the grading is
$\deg (x\la^{k}\otimes 1)=\deg (1\otimes x\mu^{k})=k$ if $x\in
sl_2(\CC)$ and $k\in\ZZ$, and
$\deg\la=\deg\mu=1$). Indeed, by this interpretation the $r$-matrices
$r_i$, $i=1,2,3$ correspond respectively to
$$
\bar r_1={1\over 2}({1\over{\la-\mu}}-{1\over{\mu-\la}})t,
\bar r_2={1\over 2}({{\la+\mu}\over{\la-\mu}}-{{\la+\mu}\over{\mu-\la}})t
+2(e\otimes f-f\otimes e),$$
$$
\bar r_3={1\over 2}({{\la\mu}\over{\la-\mu}}-{\la\mu\over{\mu-\la}})t
+2\la e\otimes f -2\mu f\otimes e
\eqno(2)$$
(we use the convention ${1\over{\la-\mu}}=\sum_{i\ge 0}\mu^i/\la^{i+1}$).

Let us now consider the tensors
$$
t_1={1\over 2}({1\over{\la-\mu}}+{1\over{\mu-\la}})t,
t_2={1\over 2}({{\la+\mu}\over{\la-\mu}}+{{\la+\mu}\over{\mu-\la}})t,
t_3={1\over 2}({{\la\mu}\over{\la-\mu}}+{\la\mu\over{\mu-\la}})t
$$

They can be considered as the split Casimir elements corresponding to
the invariant pairings on $sl_2(\CC((\la^{-1})))$
$$
\langle a(\la),b(\la)\rangle_1=\res_{\la=0}\langle a(\la), b(\la)\rangle
d\la,
\langle a(\la),b(\la)\rangle_2={1\over 2}\res_{\la=0}\langle a(\la),
b(\la)\rangle {d\la \over\la},$$
$$
\langle a(\la),b(\la)\rangle_3=\res_{\la=0}\langle a(\la),
b(\la)\rangle{d\la
\over\la^2},
$$
which allow the interpretation of $(sl_2(\CC((\la^{-1}))), \bar r_i)$ as
double algebras
(see \cite{D}); here $\langle,\rangle$ is the Killing form on $sl_2(\CC)$.

\proclaim{Lemma}
Let $a_i$ be numbers ($i=1,2,3$), then we have
$$
[\![\sum_i a_i\bar r_i,\sum_ia_i\bar r_i]\!]=[\![\sum_i
a_it_i,\sum_ia_it_i]\!].
\eqno(3)$$
\endproclaim

\demo{Proof} We first show (3) for $a_{1}=a_{2}=0$: in this case, it
follows from the structure of the Manin triple in the rational case; we
can also see it as a consequence of the identity (in formal series)
${{\la\mu}\over{\la-\mu}}{{\mu\nu}\over{\mu-\nu}}+\text{p.c.}=0$ (itself
a consequence of
${1\over{\mu^{-1}-\la^{-1}}}{1\over{\nu^{-1}-\mu^{-1}}}+\text{p.c.}=0$).
We
prove (3) for $a_{1}=a_{3}=0$ and $a_{2}=a_{3}=0$ similarly. Returning now
to the identity for $a_{2}=a_{3}=0$, the action of $\tau_{E}$
transforms
$t_3$ into $t_3+Et_2+E^2t_1$, and $\bar r_3$ into $\bar r_3+E\bar
r_2+E^2\bar r_1$,
so we have
$[\![\bar r_3+E\bar r_2+E^2\bar r_1, \bar r_3+E\bar r_2+E^2\bar r_1]\!]=
[\![t_3+Et_2+E^2t_1,
t_3+Et_2+E^2t_1]\!]$; so we deduce from
$[\![\bar r_{i},\bar r_{j}]\!]=[\![t_{i},t_{j}]\!]$ for $i,j\le 2$ or
$i,j\ge 2$, the same identity for all $i$, $j$. \qed
\enddemo

\subheading{4. Manin triples}

Let us now consider $g=sl_2(\CC((\la^{-1})))$, endowed with the Poisson-Lie
structure given by $\sum a_ir_i$. Our aim in this section is to describe
it as a double algebra.
For this, and following \cite{STS}, we study the eigenvalues of the
operator $R:g\to g$ defined
by $RA=\langle \sum_i a_ir_i, A\otimes 1 \rangle_{(a_i)}$, where
$\langle,\rangle_{(a_i)}$ is the scalar product corresponding to $\sum
a_it_i$

Let us compute $\langle,\rangle_{(a_i)}$. $\sum_i a_it_i=\sum a'_{i+j+2}
(2e\la^i\otimes f\mu^j+2f\la^i\otimes e\la^j+h\la^i\otimes h\mu^j)$, with
$a'_{1,3}=a_{1,3}$ and $a'_{2}=2a_{2}$. The inverse of the matrix $A$ with
coefficients $A_{ij}=a_{i+j+2}$ has coefficients $(A^{-1})_{ij}=
b_{i+j}$, where $b_n$ is the coefficient of $\la^{-n}$ in the expansion of
$(a_1\la^{-1}+a_2+a_3\la)^{-1}$. We will have
$\langle e\la^i,f\mu^j\rangle_{(a_i)}
={1\over 2}\langle h\la^i,h\mu^j\rangle_{(a_i)}
=b_{i+j}$, and so
$$
\langle A(\la),B(\la)\rangle_{(a_i)}=\res_{\la=\infty}
\langle A(\la),B(\la)\rangle
{d\la\over{a_1+2a_2\la+a_3\la^2}}\eqno(4)
$$
We can then express $R$ as
$$\eqalign{
R(A)(\mu) &=\langle r(\la,\mu),A(\la) \rangle_{(a_i)}\cr
&=\res_{\la=\infty}\{
t({a_1+(\la+\mu)a_2+\la\mu a_3}){1\over 2}({1\over{\la-\mu}}-{1\over{\mu-\la}})
+2(a_2+a_3\la)e\otimes f\cr
&-2(a_2+a_3\mu)f\otimes e
\}A(\la){d\la\over{a_1+2a_2\la+a_3\la^2}}
\cr
&=\res_{\la=\infty}tA(\la){1\over 2}({1\over{\la-\mu}}-{1\over{\mu-\la}})d\la
+\res_{\la=\infty}{{-t(a_2+a_3\la)}\over{a_1+2a_2\la+a_3\la^2}}A(\la)d\la\cr
&+\res_{\la=\infty}
{{2e\otimes f(a_2+a_3\la)}\over{a_1+2a_2\la+a_3\la^2}}A(\la)d\la
-\res_{\la=\infty}
{{2f\otimes e(a_2+a_3\mu)}\over{a_1+2a_2\la+a_3\la^2}}A(\la)d\la,
\cr}$$
where $\langle,\rangle$ is understood in the two last lines. Since
$\langle t,A\rangle=2A$ ($\langle a,b\rangle$ is the trace of $ab$ in the
fundamental representation), and denoting for a formal series $\phi=\sum
\phi_i\la^i$, $\phi_{>0}=\sum_{i>0} \phi_i\la^i$ and
$\phi_{\le 0}=\sum_{i\le 0} \phi_i\la^i$,
we get for the first term $A(\mu)_{>0}-A(\mu)_{\le0}$.

We see that $R$ acts as $1$ on $a\otimes\CC[[\la^{-1}]]$, and
$-1$ on $f\otimes \la \CC[\la]$. In addition, the subspace
$h\otimes \CC[\la]\oplus e\otimes \la^{-1}\CC[\la]$ is stable under $R$.
Let $A=\sum_{i\ge 0}\a_i h\la^i+\sum_{i\ge -1}\b_ie\la^i$ belong to this space.
We now restrict ourselves to the particular case $a_1=1$, $a_2=0$ and $a_3=-1$.
Pose $RA=\sum_{i\le 0}\a'_i h\la^i+\sum_{i\le 1}\b'_ie\la^i$.
Then
$$
\a'_0=\a_0+2(\a_{2}+\a_{4}+\cdots), \b'_0=\b_0+2(\b_{2}+\b_{4}+\cdots),
\b'_1=\b_{-1}+2(\b_{1}+\b_{3}+\cdots),
$$
and
$$\a'_i=-\a_i, \quad \b'_i=-\b_i \text{ for } i<0.
$$
The space is then the sum of the eigenspaces $\CC h\oplus\CC e\oplus\CC
e\la^{-1}$ for the
eigenvalue $1$, and $\{A(\la)|\a_0+\a_{2}+\cdots=0,\b_0+\b_{2}+\cdots=0,
\b_{-1}+\b_{1}+\b_{3}\cdots=0\}$ for $-1$.
So we have:

\proclaim{Proposition} The eigenspaces of $R$ are $a[[\la^{-1}]]$ and $g_+$,
where $g_+$ is the
subspace of $a\otimes \CC(\la)$
of functions $A(\la)$ with only
poles zero and infinity, and such that $A(1)$ and $A(-1)\in b_-$,
$[A(1)+A(-1)]_h=0$, and $A\in{1\over \la} n_+\oplus
b_+\oplus\la a[[\la]]$
(in the
completion at $0$ $a\otimes \CC((\la))$). Here $n_+=\CC
e$, $\underline h=\CC h$, $b_+=\underline h\oplus n_+$,
$b_-=\CC h+\CC f$, $[\  ]_h$ is the natural projection from $b_-$ to
$\underline h$.

These two spaces are supplementary Lie subalgebras
of $g$, isotropic
for the form
$$\langle A(\la),B(\la) \rangle_{(1,0,-1)}=\res_{\la=\infty}\langle
A(\la),B(\la)\rangle {d\la\over{1-\la^2}}.
$$
So the decomposition $g=g_{+}\oplus a[\la]$, together with the
above scalar product on $g$, defines a Manin triple.
\endproclaim

These last facts can be checked directly, for example for the isotropy of
$g_+$, the residue formula gives:
if $A$ and $B\in g_+$, then
$$\eqalign{
\langle A(\la),B(\la)\rangle_{(a_i)}&= -\res_{\la=1}\langle A,B\rangle
{d\la\over{\la-1}}{{-1}\over{\la+1}}-\res_{\la=-1}\langle A,B\rangle
{d\la\over{\la+1}}{-1\over{\la-1}}\cr
&-\res_{\la=0}\langle A, B\rangle
{d\la\over{1-\la^2}}\cr
&={1\over 2 }(\res_{\la=1}\langle A,B\rangle{d\la\over{\la-1}}-
\res_{\la=-1}\langle A,B\rangle{d\la\over{\la+1}})\cr
&-\res_{\la=0}
\langle{1\over\la}n_+^{A}+b_+^{A}+\cdots, {1\over\la}n_+^{B}+b_+^{B}+\cdots
\rangle d\la(1+\la^2+\cdots)}
$$
The first term is zero by the compatibility conditions at points $1$ and $-1$,
and the second vanishes also.

\demo{Remarks}

1. Remind, that a solution $r(\la,\mu)$ of CYBE is called non-degenerate if it
is non-degenerate for at least one point $(\la,\mu)$. Drinfeld conjectered that
every non-degenerate rational $r$-matrix $r(\la,\mu)$ is equivalent
to a solution of CYBE of the form
$$
r(\la,\mu) = {t\over{\la-\mu}} + P(\la,\mu),
$$
where $P(\la,\mu)$ is a polynomial of degree at most one in each variables.
This problem was investigated by A. Stolin in \cite{S} and it follows from
his results that in the $sl_{2}(\CC)$-case there are two rational
$r$-matrices:
$$
r_{1}'(\la,\mu) = {t\over{\la-\mu}} + 2(h\otimes e - e\otimes h), \eqno(5)
$$
$$
r_{1}''(\la,\mu) = {t\over{\la-\mu}} + 2(\mu h\otimes e - \la e\otimes h).
\eqno(6)
$$
We can perform the manipulations of sect. 1 on the above matrices. Starting
with the $r$-matrices (5) and (6), we obtain the $r$-matrices
$$
r_{3}^{'}(\la,\mu) = {\la\mu t\over{\la-\mu}} + 2(\la e\otimes f + \mu
f\otimes e) + 2(\mu h\otimes e -\la e\otimes h)
$$
and
$$
r_{3}^{''}(\la,\mu) = {\la\mu t\over{\la-\mu}} +
2(\la e\otimes f + \mu f\otimes e) +2(h\otimes e - e\otimes h),
$$
which can be included in the similar compatible triple of multiplicative
Poisson structures (the bialgebra structure corresponding to $r_{1}'$
has been quantized recently in [KhST]).

2. In the case of different parameters $(a_{i})$, the points $1$ and $-1$
should be shifted to some other points $z_i$
of $\CC^*$ (in the same time the form ${{d\la}\over{\la^2-1}}$ is
changed to ${{d\la}\over{(\la-z_1)(\la-z_2)}}$). It would be interesting
to see how this family a Manin triples degenerates to the known cases
rational and trigonometric ones (as they are described in \cite{D}).

3. It could also be interesting to generalize this to a situation
with more marked
points. We see that the description given here does not fit in the
adelic framework if \cite{D} (the form is not the sum of residues at all
singular points
but only the residue at one point, and the ``evident'' half of Manin pair is
not the algebra of functions regular outside these points).
\enddemo


\subheading{Acknowledgements} The authors express their gratitude to
V. Fock for discussions related to the matter of this text, and to
A.Yu. Orlov for discussions concerning the Poisson brackets of
monodromy operators in the higher Hamiltonian structures.

This paper was partly written during the
visit paid by the second author to the Ecole Polytechnique, to whom he
would like to express his thanks.
The second author is supported by ISF grants
MGK000, and  RFFI 94-02-03379 and 94-02-14365, and also
supported by the CNRS.

\enddocument